\documentstyle[12pt]{article}
\def\Bbb{\bf}    

\def\cz{{\unitlength1pt\begin{picture}(9,8)
         \put(0,0){{\rm C}}
         \put(3,0.5){\line(0,1){7.5}}
         \end{picture}}}

\def\gz{{\Bbb Z}}
\def\pr{{\rm I \! P}}
\def\qz{{\Bbb Q}}
\newcommand{\skp}{\hspace{1pt}}

\newcommand{\ph}{\varphi}

\newcommand{\id}{{\mbox{\rm Id}\skp}}

\newcommand{\img}{{\mbox{\rm im}\skp\skp}}

\newcommand{\In}{{\mbox{\rm In}\skp}}

\newcommand{\End}{{\mbox{\rm End}\skp\skp}}
\newcommand{\Pic}{{\mbox{\rm Pic}\skp}}

\newcommand{\calc}{{\cal C}}

\newcommand{\cali}{{\cal I}}

\newcommand{\calm}{{\cal M}}
\newcommand{\caln}{{\cal N}}
\newcommand{\calo}{{\cal O}}

\newcommand{\calu}{{\cal U}}

\newtheorem{prop}{Proposition}[section]
\newtheorem{theorem}[prop]{Theorem}
\newtheorem{lemma}[prop]{Lemma}

\newtheorem{defi}[prop]{Definition}
\newcommand{\pf}{{\em Proof. }}
\newcommand{\qed}{{\hfill$\Diamond$}\vspace{1.5ex}}

\newcommand{\modg}{{\calm_g}}
\newcommand{\ch}{{\mbox{\rm ch}\skp}}
\newcommand{\arctanh}{{\mbox{\rm arctanh}}}

\title{Recursive relations for the cohomology ring of moduli
spaces of stable bundles}
\author{Bernd Siebert\\
        Mathematisches Institut\\
        Bunsenstr.~3--5\\
        D-37073 G\"ottingen\\[5pt]
        siebert@cfgauss.uni-math.gwdg.de
        \and
        Gang Tian\\
        Courant Institute\\
        251 Mercer Street\\
        New York, NY 10012, USA\\
        tiang@cims.nyu.edu}
\begin{document}
\maketitle


\noindent
The moduli spaces under consideration basically arise in two different
ways:  Algebraically as a space of (isomorphism classes) of
representations $\rho:\pi_1(S_g)\rightarrow\pr U(2)\simeq SO(3)$ which
{\em do not lift to $SU(2)$}. Here $\pi_1(S_g)$ is the fundamental
group of the closed oriented surface $S_g$ of genus $g$ ($\ge2$ in the
following). Or geometrically as space of (isomorphism classes) of
stable bundles $E$ over a Riemann surface $\Sigma$ of genus $g$ with
fixed determinant bundle $L=\det E$  of {\em odd} degree. It was an
observation of Mumford \cite{mumford} that in the latter instance the
moduli space is a non-singular projective variety (of dimension
$3g-3$). Let us denote this variety by ${\calm}(\Sigma,L)$ to indicate
its dependence on the particular choices of $\Sigma$ and $L$. The space
of representations on the other hand depends only on $g$ and has the
structure of a differentiable manifold which we denote by ${\caln}_g$
\cite{nase}. Then the celebrated theorem of Narasimhan and Seshadri
says that ${\caln}_g$ is the manifold underlying ${\calm}(\Sigma,L)$
for any $\Sigma$, $L$ \cite{nase}.

We are dealing here with the simplest non-trivial moduli spaces of
stable bundles (they are non-singular, complete, Fano, with Picard
group infinite cyclic of rank 1) and quite a bit is known about their
cohomology.  For example there are several ways of computing the
Betti-numbers \cite{newstead1}, \cite{harder}, \cite{atiyah-bott},
\cite{kirwan0}, the cohomology is known to be torsion-free
\cite{atiyah-bott} and explicit generators for the integral ring have
been found \cite{newstead1}, \cite{atiyah-bott}. More recently Thaddeus
\cite{thaddeus} found a beautiful formula for the highest intersection
products relating to the Verlinde formula. Around the same time Kirwan
succeeded in proving Mumford's conjecture that a certain (huge)
canonical set of relations for the cohomology ring is complete.

What is unsatisfactory about the picture is that still a big
computational effort is required to find generators for the relations
explicitely. And already low genus, say $g=5$, demands the use of
computers which in turn will not reach beyond $g=15$ or so.

The content of this paper remedies this situation and should be good
news for everybody presently using ${\caln}_g$: The basic part of
$H^*({\caln}_g)$ (the subring generated by Newstead's classes $\alpha$,
$\beta$, $\gamma$) is a complete intersection ring and there is a very
simple inductive formula for generators of the relation ideal! Namely,
\begin{theorem}
$\langle\alpha,\beta,\gamma\rangle=\qz[\alpha,\beta,\gamma]/
(f_1^g,f_2^g,f_3^g)$ with
\[
  f_1^{g+1}\ =\ \alpha f_1^g+g^2 f_2^g,\ \
  f_2^{g+1}\ =\ \beta f_1^g+\frac{2g}{g+1}f_3^g, \ \
  f_3^{g+1}\ =\ \gamma f_1^g
\]
for $g\ge1$ and formally setting
$(f_1^1,f_2^1,f_3^1)=(\alpha,\beta,\gamma)$.  The $f_i^g$ are uniquely
determined by their respective initial terms $\alpha^g$,
$\alpha^{g-1}\beta$ and $\alpha^{g-1}\gamma$.
\end{theorem}
As a consequence the monomials $\alpha^a\beta^b\gamma^c$ with $a+b+c<g$
form a $\qz$-basis for this subring. Moreover, there is a generating
function $\Phi(t)$ for generators of {\em all} the relation ideals,
i.e.\ simultaneously for all $g$,
s.th.\ $(f_1^g,f_2^g,f_3^g)=(\Phi^{(g)},\Phi^{(g+1)},\Phi^{(g+2)})$ as
ideals ($\Phi^{(r)}$ the $r$-th derivative at $t=0$). $\Phi$ is
characterized by the simple (formal) functional equation
\[
  \Phi'(t)\ =\ \frac{\alpha+\beta t+2\gamma t^2}{1-\beta t^2}\cdot\Phi(t)\,.
\]
For an explicit expression of $\Phi$, cf.\ Definition~\ref{genfct}.

Besides its aesthetical appeal, why is this exciting? First of all our
approach is technically simpler and yields even more insight into the
cohomology ring than any other approach: The only input we have to use
is the embedding $\ph:{\calm}(\Sigma,L)\rightarrow G(g+3,2g+2)$ of
Desale and Ramanan ($\Sigma$ hyperelliptic) together with an intrinsic
characterization of the pull-back of the tautological bundle
\cite{ramanan}; and the total dimension of $H^{2*}({\caln}_g)$ for
which various methods of computation are available. From the recursion
relations one should actually be able to derive without much additional
effort more or less everything that is known on the cohomology of
${\caln}_g$, e.g.\ the vanishing of the Pontryagin ring and the Chern
classes of ${\caln}_g$ above degree $4g-4$ \cite{kirwan},
\cite{thaddeus}, \cite{gieseker}.

Our method should also be applicable to moduli spaces of bundles of
even degree, once the technical problem of non-existence of a universal
bundle in the literal sense is overcome \cite{ramanan0}.  Furthermore,
there are chances that the method of proving completeness of certain
relations by an explicit computation of the Leitideal generalizes
to moduli spaces of higher rank bundles (where none of the other
methods apply at the moment).

Recently the cohomology of $\caln_g$ occurred as instanton Floer
homology on the three-manifold $S^1\times S_g$ \cite{salamon}. It might
well be that there are applications of our formulas to this circle of
ideas or even to Donaldson invariants via gluing formulas
\cite{donaldson}.

Finally, and this was the main motivation for the authors to search for
a minimal set of relations, finding a good presentation for the
cohomology ring is the first essential step in the computation of the
quantum cohomology of ${\caln}_g$, cf.\ \cite{siebert-tian}. The latter
should compute the so-called Fukaya-Floer homology of a trivial product
$S^1\times S_g$, cf.\ \cite{donaldson2}. Classically there follows at
least a new integral formula for the highest intersection products on
${\caln}_g$, hence for Bernoulli numbers, by the formula of
\cite{thaddeus}, cf.\ \cite{siebert-tian}.

After this paper had been finished and put on the electronic bulletin
board {\tt alg-geom} the authors learnt of work of several other
mathematicians proving similar results \cite{zagier},
\cite{baranovsky}, \cite{king}. There preprints are all not available
at the moment, but we can say the following: V.Baranovsky also uses the
Desale-Ramanan model to get relations, while D.Zagier starts with the
original Mumford relations of minimal degree. V.Balaji, A.King and
P.E.Newstead in turn set up the recursion relations in their form from
a geometric argument and compute the coefficients by studying explicit
families of bundles. Previously, P.E.Newstead had computed the relations
modulo $\beta$ \cite{newstead3}.
\vspace{5ex}

\noindent
{\bf Acknowledgements:} We thank A.King for kindly informing us of
\cite{zagier}, \cite{baranovsky} and \cite{king}. The first named
author wants to thank the Deutsche Forschungsgemeinschaft for their
support during the academic year 1993/94 and the Courant Institute
(where this work was done) for hospitality. The second named author
is partly supported by NSF grants and an Alfred P.\ Sloan fellowship.

\section{Method and notation}
Let us now fix a Riemann surface $\Sigma$ of genus $g$ and a line
bundle $L$ on $\Sigma$ of odd degree. We write
$\modg={\calm}_g(\Sigma,L)$.  Generators for $H^*(\modg,\qz)$ occur as
coefficients in the K\"unneth decomposition of a characteristic class
associated to the universal bundle $\calu$ over $\modg\times\Sigma$
\[
  c_2(\End\calu)\ =\ -\beta+4\psi+2\alpha\otimes\omega\,.
\]
Here $\omega\in H^2(\Sigma,\gz)$ is the normalized volume form and $\psi$
is the part of type $(3,1)$. Choosing a standard basis $e_i$ for
$H^1(\Sigma,\gz)$, $i=1,\ldots,2g$ (s.th.\  $e_i e_j=0$ unless $i\equiv
j(g)$ and $e_i e_{i+g}[\Sigma]=1$), $\psi$ decomposes further
$\psi=\sum_{i=1}^{2g} \psi_i\otimes e_i$. The classes $\alpha\in
H^2(\modg,\qz)$, $\psi_i \in H^3(\modg,\qz)$ and $\beta\in
H^4(\modg,\qz)$ are actually integral and generate the cohomology ring
\cite{newstead1}, \cite{atiyah-bott}.

There is an interesting subring of $H^*({\caln}_g,\qz)$ to which the
intersection pairing may easily be reduced by geometric arguments as
noted by Thaddeus \cite{thaddeus}. Namely, any orientation preserving
diffeomorphism of $\Sigma$ induces a diffeomorphism of ${\caln}_g$ by
acting on $\pi_1$. The corresponding action on $H^*(\modg,\qz)$ leaves
$c_2(\End\calu)$ fixed. Thus $\alpha$ and $\beta$ are invariant and the
$\psi_i$ transform dually to the $e_i$. It is not hard to show that the
ring of such transformations is precisely the subring generated by
$\alpha$, $\beta$ and a newly defined class
$\gamma:=2\sum_{i=1}^{2g}\psi_i\psi_{i+g}$, or more intrinsically
$\psi^2=\gamma\omega$ (in view of the functional equation for the
generating function $\Phi$ it might be more natural to take twice this
class, but to be consistent with the work of Newstead and Thaddeus we
keep this definition).

On the other hand, there is a method introduced by Mumford to construct
relations among the generators, cf.\ \cite{atiyah-bott}: Letting $L$
vary among the line bundles of fixed (odd) degree $D$ (the space of
which we denote by $\Pic^D(\Sigma)$), one gets a moduli space
$\tilde{\calm}(\Sigma)$ and a fibration
$\tilde{\calm}(\Sigma)\rightarrow\Pic^D(\Sigma)$ with fibre $\caln_g$.
In rational cohomology this fibration is trivial, i.e.\
$H^*(\tilde{\calm}(\Sigma),\qz)\simeq H^*(\caln_g,\qz)\otimes_\qz
H^*(\Pic^D(\Sigma),\qz)$. But $\Pic^D(\Sigma)\simeq\Pic^0(\Sigma)$ and
$H^*(\Pic^0(\Sigma),\qz)$ is a free alternating algebra in $2g$
generators ${\hat e}_i$ of degree 1 (the push-forwards of $e_i$ under
the Jacobi map $\Sigma\rightarrow\Pic^0(\Sigma)$). Now let
$\tilde\calu$ be the universal bundle over
${\tilde\calm}_g\times\Sigma$ and let $\tilde\pi:{\tilde\calm}_g
\times\Sigma\rightarrow{\tilde\calm}_g$ be the projection. If $D=4g-3$,
$R^1\tilde\pi_*\tilde\calu=0$ and $\tilde\pi_*\tilde\calu$ is locally
free of rank $2g-1$. By Grothendieck-Riemann-Roch then the Chern
classes of $\tilde\pi_*\tilde\calu=\tilde\pi_!\,\tilde\calu$ are
expressed as polynomials in $\alpha$, $\beta$, $\psi_i$ and $\hat
e_i$.  Equating to zero the coefficients of $\hat e_{i_1}\ldots\hat
e_{i_\nu}$ in $c_r(\tilde\pi_*\tilde\calu)$, $r\ge 2g$ (which vanish
for rank reasons) gives a number of relations among the generators
$\alpha$, $\beta$ and $\psi_i$. The point of letting $L$ vary is of
course to lower the degrees of the relations by up to $2g$. The
smallest degree of a relation we thus obtain is $4g-2g=2g$.

The conjecture of Mumford which Kirwan recently succeeded to prove as
remarked in the introduction, is that this set of relations is
complete.  But in view of the explicit formula for the intersection
pairing of Thaddeus the authors could not believe in this end of the
story. In fact, computer evidence (up to $g=18$ using ``Macaulay''
\cite{mcly}) showed that for the subring generated by $\alpha$,
$\beta$, $\gamma$ there should be only three independent relations of
degrees $g$, $g+1$, $g+2$ respectively, coming from the lowest degree
equation $c_{2g}=0$. Unfortunately, explicit computations are rather
arduous, e.g.\ due to the presence of the odd degree classes $\hat
e_i$.

To find relations of low degrees without bringing $\Pic(\Sigma)$ into
the game we first remark that because the ${\calm}(\Sigma,L)$ are all
diffeomorphic as long as we fix the genus of $\Sigma$ and the degree of
$L$ modulo 2, we may restrict ourselves to a hyperelliptic curve
$\Sigma$. In this case, there is a closed embedding
$\ph:\modg\hookrightarrow G(g+3,2g+2)$ into a Grassmannian as follows
\cite{ramanan}: Let $p:\Sigma\rightarrow\pr^1$ represent $\Sigma$ as
two-fold covering of $\pr^1$, $\iota:\Sigma\rightarrow\Sigma$ the
corresponding hyperelliptic involution (s.th.\ $p\circ\iota=p$) and
$B\subset\pr^1$ the branch locus of $p$ ($\sharp B=2g+2$). Now let $L$
be a line bundle of degree $2g+1$ (as opposed to $d=4g-3$ in Mumford's
method) and let $E$ be a stable 2-bundle over $\Sigma$ with determinant
$L$. Applying $\iota^*$ to $E\otimes\iota^* E$ and switching factors
induces an involution $J:  E\otimes\iota^*E\rightarrow E\otimes\iota^*
E$. Denote by $p_*(E\otimes \iota^*E)^\natural$ the $J$-anti-invariant
subsheaf of $p_*(E\otimes\iota^* E)$.  One shows
$h^0(\pr^1,p_*(E\otimes\iota^* E)^\natural)=g+3$ (loc.\ cit.,
Prop.~2.2). In a branch point $t\in B$ we have the identification
\[
  p_*(E\otimes\iota^* E)^\natural_t\ \simeq\ p_*(E\wedge E)_t
  \ =\ (p_*L)_t\,.
\]
The map $\ph:\modg\rightarrow G(g+3,2g+2)$ is then defined as
\[
  E\longmapsto\left(H^0\Big(\pr^1,p_*(E\otimes\iota^* E)^\natural
  \Big)\Big|_B\subset H^0(\pr^1,p_* L|_B)\simeq\cz^{2g+2}\right)\,.
\]
Now let $S$ and $Q$ be the universal bundle and the universal quotient
bundle on $G(g+3,2g+2)$ respectively. The key observation is that $Q$
has rank $g-1$!  Similarly to Mumford's method we ``just'' have to
express the Chern classes of $\ph^* Q$ in terms of $\alpha$, $\beta$,
$\gamma$ (essentially by Grothendieck-Riemann-Roch of course) to get
relations $c_r(\ph^* Q)=0$, $r=g,g+1,g+2$.

The Chern class computations are based on the following exact
sequence.
\begin{lemma}\label{exact}
  On $\modg\times\Sigma$, there is an exact sequence
  \[
    0\longrightarrow\hat p^*\hat p_*(\calu\otimes\iota^*\calu)^\natural
    \longrightarrow\calu\otimes\iota^*\calu\longrightarrow
    S^2\calu|_{\modg\times p^{-1}(B)}\longrightarrow 0\,,
  \]
  where $\hat p=\id\times p:\modg\times\Sigma\rightarrow\modg\times\pr^1$.
\end{lemma}
\pf
The restriction map $\hat p^*\hat p_* E\rightarrow E$ is injective for
any locally free sheaf $E$ by flatness of $p$. Composing this map with
the inclusion $\hat p^*\hat
p_*(\calu\otimes\iota^*\calu)^\natural\hookrightarrow \hat p^*\hat
p_*(\calu\otimes\iota^*\calu)$ yields exactness at the first place.
Outside of the branch points this map is obviously an isomorphism (by
anti-invariance elements of $\hat p^*\hat p_*(\calu\otimes
\iota^*\calu)^\natural$ are uniquely determined by their behaviour on
one branch). It is then a matter of linear algebra to check that at a
branch point $y\in\Sigma$ the cokernel is given by mapping a germ of
sections of $\calu \otimes\iota^*\calu$ to
$(s+J(s))(y)\in\calu\otimes\iota^*\calu |_{\modg\times p^{-1}(B)}\simeq
S^2\calu|_{\modg\times p^{-1}(B)}$.
\qed

\section{Computations of Chern classes}
This section contains the computational heart of the paper, summarized
in the following proposition. We adopt the convention that in writing
Chern classes or Chern characters as analytic functions in certain
cohomology classes we understand to evaluate on these classes the
corresponding power series expansion.
\begin{prop}\label{chernquot}
  Denote by $c(\ph^*Q)=\sum_{i\ge0}c_i(\ph^*Q)$ the total Chern class
  of $\ph^*Q$. Then
  \[
    c(\ph^*Q)=(1-\beta)^{-1/2}\exp\bigg[\alpha+\Big(\alpha+
    \frac{2\gamma}{\beta}\Big)\sum_{m\ge1}\frac{\beta^m}{2m+1}\bigg]\,.
  \]
\end{prop}
Note that the $c_i(\ph^*(Q))$ are really {\em polynomials} in $\alpha$,
$\beta$, $\gamma$, the denominator $\beta$ cancels.  Also, if one
prefers, one could write $(\arctanh(\sqrt{\beta})/\sqrt{\beta})-1$
(with $\sqrt{\beta}$ formally adjoint to $H^*(\modg)$) instead of the
infinite sum.

Before turning to the proof we will need some preparations.  Letting
$\pi:\modg\times\pr^1\rightarrow\modg$ and $\tilde\pi=\pi\circ\hat
p:\modg\times \Sigma\rightarrow\modg$ be the projections, we know
$\ph^*S={\tilde\pi}_*( \calu\otimes\iota^*\calu)^\natural\otimes
H^{-1}$ with $H\in\Pic(\modg) \simeq\gz$ the ample generator and
${\tilde\pi}_*( \calu\otimes\iota^*\calu)^\natural:=\pi_*(\hat p_*
(\calu\otimes\iota^*\calu)^\natural)$ \cite{ramanan}. To apply
Grothendieck-Riemann-Roch to this sheaf we first need a closed formula
for the Chern character of $\calu\otimes\iota^*\calu$. We will use a
compuational trick which the authors learned from \cite[...]{kirwan} to
represent most of the terms in exponential form. This will be
convenient later when we transform back to Chern classes. For that fix
a large number $N$ such that $\beta^N=0$ (e.g.\ for dimension reasons)
and let $\mu_1,\ldots,\mu_N\in\cz$ be such that
\[
  N_k(\mu_1,\ldots,\mu_N)\ :=\ \sum_{\nu=1}^N\mu_\nu^k\ =\ \frac{1}{k+1}
  \ \ \ \mbox{for } 1\le k\le N.
\]
The existence of $\mu_1,\ldots,\mu_N$ is clear either by an elementary
argument or from the finiteness of the map
$(N_1,\ldots,N_N):\cz^N\rightarrow\cz^N$.
\begin{lemma}
Formally adjoining $\sqrt{\beta}$ and $\alpha':=\alpha+2\gamma/\beta$
(both of real degree 2) to $H^{2*}(\modg,\qz)$ the following holds:
\begin{eqnarray*}
  \ch(\calu\otimes\iota^*\calu)&=&e^\alpha\bigg[\Big(
  2+e^{\sqrt{\beta}}+e^{-\sqrt{\beta}}\Big)(1+D\omega)-2\alpha\omega\\
  &&-\,\alpha'\omega\sum_\nu
  \Big(e^{\mu_\nu\sqrt{\beta}}+e^{-\mu_\nu\sqrt{\beta}}-2\Big)\bigg]\,.
\end{eqnarray*}
\end{lemma}
\pf
We have
$c(\calu^*\otimes\calu)=1+c_2(\End\calu)=1-\beta+4\psi+2\alpha\omega$
with $\psi=\sum_i\psi_i\otimes e_i$. Note that
$\iota^*:H^1(\Sigma)\rightarrow H^1(\Sigma)$ is just multiplication by
$-1$. In fact, if $\delta$ is a closed $1$-form on $\Sigma$, then
$\delta+\iota^*\delta$ is closed and $\iota^*$-invariant, hence
$\delta+\iota^*\delta=p^*df=dp^*f$ for $f\in\calc^\infty(\pr^1)$ since
$H^1(\pr^1)=0$. But any orientation preserving diffeomorphism leaves
$c(\End\calu)$ unchanged, hence $\iota^*\psi_i=-\psi_i$.  Now it is
almost clear and can be easily verified by a standard Chern class
computation that
$c_2(\calu^*\otimes\iota^*\calu)=-\beta+2\alpha\omega$, i.e.\ that the
factor $4\psi$ drops out. Instead, there is a non-trivial $c_4$, namely
$c_4(\calu^*\otimes\iota^*\calu)=4\psi^2=4\gamma\omega$. Summarizing,
we have
\[
  c(\calu^*\otimes\iota^*\calu)\ =\ 1+(-\beta+2\alpha\omega)
  +4\gamma\omega\,.
\]
Next, for any bundle $E$ with only $c_2$ and $c_4$ non-vanishing
\[
  \ch_{2k}(E)\ =\ 2(-1)^k c_2(E)^{k-2}[c_2(E)^2-kc_4(E)],\ \ \
  \ch_{2k+1}(E)=0
\]
(induction on $k$). Thus for $k\ge2$ ($k=1$: $\ch_2(\calu^*\otimes
\iota^*\calu )=-2c_2=2\beta-4\alpha\omega$)
\begin{eqnarray*}
  \ch_{2k}(\calu^*\otimes\iota^*\calu)&=&2(-1)^k\left((-\beta)^{k-2}
  +(k-2)(-\beta)^{k-3}2\alpha\omega\right)(\beta^2-4\alpha\omega-4k
  \gamma\omega)\\
  &=&2\beta^k-4k(\alpha\beta+2\gamma)\beta^{k-2}\omega\,.
\end{eqnarray*}
Formally adjoining $\sqrt{\beta}$ and $\alpha'=\alpha+2\gamma/\beta$
($\beta$ in the denominator always cancels in the following) we get
\begin{eqnarray*}
  \ch(\calu^*\otimes\iota^*\calu)&=&2+2\sum_{k\ge0}\frac{\beta^k}{(2k)!}
  -2(\alpha\beta+2\gamma)\omega\sum_{k\ge2}2k\frac{\beta^{k-2}}{(2k)!}
  -\frac{4\alpha\omega}{2}\\
  &=&2+e^{\sqrt{\beta}}+e^{-\sqrt{\beta}}-2\alpha\omega-2\alpha'\omega
  \sum_{k\ge1}\frac{1}{2k+1}\frac{\beta^k}{(2k)!}\,.
\end{eqnarray*}
The computational trick consists in writing
\[
  2\sum_{k\ge1}\frac{1}{2k+1}\frac{\beta^k}{(2k)!}\ =\
  \sum_{\nu=1}^N\Big(e^{\mu_\nu\sqrt{\beta}}
  +e^{-\mu_\nu\sqrt{\beta}}-2\Big)\,.
\]
Finally using $c_1(\calu)=\alpha+D\omega$ ($D=2g+1$) together with the
isomorphim $\calu\simeq\calu^*\otimes\det\calu$ and the
multiplicativity of the Chern character we deduce
$\ch(\calu\otimes\iota^*\calu)=e^\alpha
\ch(\calu^*\otimes\iota^*\calu)$ which upon inserting the previous
computations gives the stated formula.
\qed
\vspace{10pt}

Pushing-forward the exact sequence from Lemma~\ref{exact} we get
\[
  0\longrightarrow \hat p_*(\calu\otimes\iota^*\calu)^\natural
  \otimes\hat p_*\calo\longrightarrow \hat p_*(\calu\otimes\iota^*\calu)
  \longrightarrow\hat p_*\left(S^2\calu|_{\modg\times p^{-1}(B)}\right)
  \longrightarrow0
\]
(for the first term apply the projection formula).
\begin{lemma}
  Denoting $\bar\omega$ the normalized volume form on $\pr^1$ the
  following holds
  \begin{eqnarray*}
    \ch\left(\hat p_*(\calu\otimes\iota^*\calu)^\natural\right)
    &=&e^\alpha\bigg[(2+e^{\sqrt{\beta}}+e^{-\sqrt{\beta}})
    (1-\frac{\bar\omega}{2})\\
    &&+\,\bar\omega\Big(-\alpha-\frac{\alpha'}{2}\sum_\nu\big(
    e^{\mu_\nu\sqrt{\beta}}+e^{-\mu_\nu\sqrt{\beta}}-2\big)
    +(g+1)\Big)\bigg]\,.
  \end{eqnarray*}
\end{lemma}
\pf
{}From the above exact sequence, we see
\[
  \ch\left(\hat p_*(\calu\otimes\iota^*\calu)^\natural\right)
  \ =\ \left(\ch \hat p_*(\calu\otimes\iota^*\calu)-\ch
  \hat p_*(S^2\calu|_{\modg\times p^{-1}(B)})\right)/\ch(\hat p_*\calo)
\]
which by Grothendieck-Riemann-Roch applied to $\hat p$ and writing
$\bar g=g-1$ and $\ch(\calu\otimes\iota^*\calu)-\ch(S^2
\calu|_{\modg\times p^{-1}(B)})=A+B\omega$, equals
\begin{eqnarray*}
  \hat p_*\left[(A+B\omega)(1-\bar g\omega)\right]/
  \hat p_*(1-\bar g\omega)
  \ =\ \hat p_*\Big(A+(B-\bar g A)\omega\Big)/(2-\bar g\bar\omega)\\
  =\ \Big(2A+(B-\bar g A)\bar\omega\Big)\frac{1}{2}(1+
  \frac{\bar g}{2}\bar\omega)\ =\ A+B\frac{\bar\omega}{2}\,.
\end{eqnarray*}
To compute the Chern character of $S^2\calu|_{\modg\times p^{-1}(B)}$
we restrict the exact sequence of Lemma~\ref{exact} to $\modg\times
p^{-1}(B)$. Then since $\hat p^*\hat
p_*(\calu\otimes\iota^*\calu)^\natural|_{\modg\times p^{-1} (B)}$
$\simeq\ \det\calu|_{\modg\times p^{-1}(B)}$ we get
\[
  \ch(S^2\calu|_{\modg\times p^{-1}(B)})\ =\ \left[\ch(\calu\otimes\iota^*
  \calu)-\ch(\det\calu)\right](2g+2)\omega\,.
\]
Thus
\begin{eqnarray*}
  A+B\omega&=&\ch(\calu\otimes\iota^*\calu)\Big(1-(2g+2)\omega\Big)
  +\ch(\det\calu)(2g+2)\omega\\
  &=&e^\alpha\bigg[(2+e^{\sqrt{\beta}}+e^{-\sqrt{\beta}})
    (1-\omega)\\
    &&+\,\omega\Big(-\alpha-\frac{\alpha'}{2}\sum_\nu\big(e^{\mu_\nu
    \sqrt{\beta}}+e^{-\mu_\nu\sqrt{\beta}}-2\big)+(g+1)\Big)\bigg]\,,
\end{eqnarray*}
hence the claim.
\qed
\vspace{10pt}

\noindent
{\em Proof of Proposition~\ref{chernquot}}: It follows from
Proposition~2.2 of \cite{ramanan} applied to $\calu\otimes\iota^*\calu$
that $R^1\pi_*(\hat p_*(\calu\otimes\iota^*\calu)^\natural)=0$. The
Grothendieck-Riemann-Roch theorem for pushing-forward
$\hat p_*(\calu\otimes\iota^*\calu)^\natural$ by
$\pi:\modg\times\pr^1\rightarrow\modg$ thus reads
\begin{eqnarray*}
  \lefteqn{\ch(\tilde\pi_*(\calu\otimes\iota^*\calu)^\natural)\ =\
  \pi_*\left(\ch(\hat p_*(\calu\otimes\iota^*\calu)^\natural)\cdot
  (1+\bar\omega)\right)}\hspace{1cm}\\
  &=&e^\alpha\bigg[\frac{1}{2}(2+e^{\sqrt{\beta}}+e^{-\sqrt{\beta}})-\alpha
  -\frac{\alpha'}{2}\sum_\nu\Big(e^{\mu_\nu\sqrt{\beta}}+
  e^{-\mu_\nu\sqrt{\beta}}-2\Big)+(g+1)\bigg]\,.
\end{eqnarray*}
Plugging in the relation expressing the pull-back of S, i.e.\
$\ph^*S=\tilde\pi_*(\calu\otimes\iota^*\calu)^\natural\otimes H^{-1}$,
we get
\begin{eqnarray*}
  \ch\ph^*Q&=&(2g+2)-\ch\ph^*S\ =\
  (2g+2)-\ch\tilde\pi_*(\calu\otimes\iota^*\calu)^\natural/\ch(H)\\
  &=&(g-1)-\frac{1}{2}\Big(e^{\sqrt{\beta}}+e^{-\sqrt{\beta}}\Big)+\alpha
  +\frac{\alpha'}{2}\sum_\nu\Big(e^{\mu_\nu\sqrt{\beta}}
  +e^{-\mu_\nu\sqrt{\beta}}-2\Big)\,.
\end{eqnarray*}
Next we need to make use of the computational trick of Kirwan again: Assume
$M\gg 0$ s.th.\ $\alpha^M=\alpha'^M=0$ and find $\lambda_\kappa\in\cz$ with
$N_1(\lambda_1,\ldots,\lambda_M)=1$, $N_k(\lambda_1,\ldots,\lambda_M)=0$
for $2\le k\le M$. Then $\alpha=\sum_\kappa(e^{\lambda_\kappa\alpha}-1)$,
$\alpha'=\sum_\kappa(e^{\lambda_\kappa\alpha'}-1)$. Inserting we get
\begin{eqnarray*}
  \ch\ph^*Q&=&(g-1)-\frac{1}{2}\Big(e^{\sqrt{\beta}}+e^{-\sqrt{\beta}}\Big)
  +\sum_\kappa(e^{\lambda_\kappa\alpha}-1)\\
  &&+\frac{1}{2}\sum_{\kappa,\nu}\left(e^{\lambda_\kappa\alpha'+\mu_\nu
  \sqrt{\beta}}+e^{\lambda_\kappa\alpha'-\mu_\nu\sqrt{\beta}}
  -e^{\mu_\nu\sqrt{\beta}}-e^{-\mu_\nu\sqrt{\beta}}-2(e^{\lambda_\kappa
  \alpha'}-1)\right)\,.
\end{eqnarray*}
This is a sum of exponentials and as such easily transformed into the
corresponding total Chern class:
\begin{eqnarray*}
  c(\ph^*Q)&=&\Big[\Big(1+\sqrt{\beta}\Big)\Big(1-\sqrt{\beta}\Big)
  \Big]^{-1/2}
  \prod_\kappa(1+\lambda_\kappa\alpha)\\
  &&\cdot\,\bigg[\prod_{\kappa,\nu}
  \frac{1+\mu_\nu\sqrt{\beta}+\lambda_\kappa\alpha'}{1+\mu_\nu\sqrt{\beta}}
  \cdot
  \frac{1-\mu_\nu\sqrt{\beta}+\lambda_\kappa\alpha'}{1-\mu_\nu\sqrt{\beta}}
  \cdot\frac{1}{(1+\lambda_\kappa\alpha')^2}\bigg]^{1/2}\!\!.
\end{eqnarray*}
To get rid of the $\lambda_\kappa$ we observe that $\sigma_k(\lambda_1,
\ldots,\lambda_M)=1/k!$ \cite[p.862]{kirwan}. The product over $\kappa$ can
thus be carried out, e.g.
\begin{eqnarray*}
  \lefteqn{\prod_{\kappa,\nu}\frac{1+\mu_\nu\sqrt{\beta}+\lambda_\kappa
  \alpha'}{1+\mu_\nu\sqrt{\beta}}\ =\
  \prod_{\kappa,\nu}\bigg(1+\lambda_\kappa\frac{\alpha'}{1+\mu_\nu\sqrt{\beta}}
  \bigg)\ =\ \prod_\nu\exp\frac{\alpha'}{1+\mu_\nu\sqrt{\beta}}}\hspace{4cm}\\
  &=&\exp\alpha'\sum_\nu\sum_{l\ge0}\Big(-\mu_\nu\sqrt{\beta}\Big)^l\ =\
  \exp\alpha'\sum_{l\ge0}\frac{(-\sqrt{\beta})^l}{l+1}\,.
\end{eqnarray*}
Inserting into our last formula we thus find (the term $(1+\lambda_\kappa
\alpha')^{-2}$ cancels the summand for $l=0$)
\begin{eqnarray*}
  c(\ph^*Q)&=&(1-\beta)^{-1/2}\exp\bigg(\alpha+\frac{\alpha'}{2}\sum_{l\ge1}
  \frac{(-\sqrt{\beta})^l+(\sqrt{\beta})^l}{l+1}\bigg)\\
  &=&(1-\beta)^{-1/2}\exp\bigg(\alpha+\alpha'\sum_{m\ge1}\frac{1}{2m+1}
  \beta^m\bigg)
\end{eqnarray*}
as claimed.
\qed
\vspace{10pt}

It is convenient to gather the relations in a generating function.
\begin{defi}\label{genfct}
  We define $\Phi\in\qz[\alpha,\beta,\gamma][\![t]\!]$ by
  \[
    \Phi(t)\ :=\ (1-\beta t^2)^{-1/2}
    \exp\bigg[\alpha t+\Big(\alpha+\frac{2\gamma}{\beta}
    \Big)t\sum_{m\ge 1}\frac{\beta^m t^{2m}}{2m+1}\bigg]\,.
  \]
\end{defi}

For later use let us also state a functional equation that $\Phi$ obeys.
This equation is actually equivalent to the recursion formula
to be proved below (Proposition~\ref{recursion}).
\begin{prop}\label{fctleqn}
  $\Phi$ obeys the following differential equation:
  \[
    \Phi'(t)\ =\ \frac{\alpha+\beta t+2\gamma t^2}{1-\beta t^2}
    \cdot\Phi(t)\,.
  \]
\end{prop}
\pf
Direct computation.
\qed
\vspace{10pt}

Let us add that with the same methods, it is not hard to deduce also
a closed formula for the Chern classes of $\caln_g$. The result is:
\[
  c(\caln_g)\ =\ (1-\beta)^g \exp\Big(\frac{-8\gamma}{1-\beta}\Big)\cdot
                 c(\ph^*Q)^2\,.
\]
Note the simple dependence on the genus!

\section{A minimal set of relations}
We have remarked in the introduction that the three generating
relations $f_1^g$, $f_2^g$, $f_3^g$ of degrees $g$, $g+1$ and $g+2$ are
uniquely determined by their initial terms $\alpha^g$,
$\alpha^{g-1}\beta$ and $\alpha^{g-1}\gamma$ respectively (w.r.t.\ the
reverse lexicographic order; we will prove this as an easy consequence
of the recursion relations, see Proposition~\ref{inf}). It is then an
exercise in calculus to find the following
\begin{defi}
  Writing $\Phi^{(r)}=\displaystyle\frac{d^r\Phi}{dt^r}(0)$
  we define for $g\ge1$
  \begin{eqnarray*}
    f_1^g&:=&\Phi^{(g)}\\
    f_2^g&:=&\frac{1}{g^2}\left(\Phi^{(g+1)}-\alpha\Phi^{(g)}\right)\\
    f_3^g&:=&\frac{1}{2g(g+1)}\left(\Phi^{(g+2)}-\alpha\Phi^{(g+1)}
             -(g+1)^2\beta\Phi^{(g)}\right)\,.
  \end{eqnarray*}
\end{defi}
\vspace{10pt}

We are now in a position to prove the recursion relations.
\begin{prop}\label{recursion}
  $(f_1^1,f_2^1,f_3^1)=(\alpha,\beta,\gamma)$ and inductively for $g\ge1$
  \begin{eqnarray*}
    f_1^{g+1}&=&\alpha f_1^g+g^2 f_2^g\\
    f_2^{g+1}&=&\beta f_1^g+\frac{2g}{g+1}f_3^g\\
    f_3^{g+1}&=&\gamma f_1^g\,.
  \end{eqnarray*}
\end{prop}
\pf
The first claim is by direct check. Next, the recursion relations for
$f_1^{g+1}$ and $f_2^{g+1}$ are immediate consequences of their
definition.  All the work is thus shifted to the innocuous looking
formula for $f_3^{g+1}$. What we have to show is the vanishing of
\begin{eqnarray*}
  \lefteqn{2(g+1)(g+2)\left(f_3^{g+1}-\gamma
  f_1^g\right)}\hspace{1cm}\\
  &=&\Phi^{(g+3)}-\alpha\Phi^{(g+2)}-(g+2)^2\beta\Phi^{(g+1)}-2(g+1)(g+2)
  \gamma\Phi^{(g)}\\
  &=&(g+2)!\left[(g+3)\ph_{g+3}-\alpha\ph_{g+2}-(g+2)\beta\ph_{g+1}
  -2\gamma\ph_g\right]
\end{eqnarray*}
with $\ph_k$ the $k$-th Taylor coefficient of $\Phi$ at $t=0$.
Multiplying with $t^{g+3}$ and taking the sum this will follow from
\[
  \sum_{g\ge1}(g+3)\ph_{g+3}t^{g+3}\ =\ \alpha
  t\sum_{g\ge1}\ph_{g+2}t^{g+2} +\beta
  t^2\sum_{g\ge1}(g+2)\ph_{g+1}t^{g+1}+2\gamma t^3\sum_{g\ge1}\ph_g t^g
\]
which is the part of order larger 3 of
\[
  t\cdot\Phi'\ =\ \alpha t\Phi+\beta t^2(\Phi\cdot t)'+2\gamma t^3\Phi
  \ =\ (\alpha t+\beta t^2+2\gamma t^3)\Phi+\beta t^3\Phi'\,,
\]
the functional equation for $\Phi$ (Proposition~\ref{fctleqn}).
\qed

\section{The Leitideal}
The decisive step in the proof of completeness of our relations is that
the Leitideal (initial ideal) can be computed completely and has a
particularly simple form. In all that follows we use the (graded)
reverse lexicographic order in $\qz[\alpha,\beta,\gamma]$ and write
$\In(f)$ ($\In(\cali)$) for the initial term of
$f\in\qz[\alpha,\beta,\gamma]$ (resp.\ the initial ideal of an ideal
$\cali\subset\qz[\alpha,\beta,\gamma]$). As warm-up let us check that
the initial terms of the $f_i^g$ are as promised in the last chapter:
\begin{prop}\label{inf}
  In the reverse lexicographic order $\In(f_1^g)=\alpha^g$,
  $\In(f_2^g)=\alpha^{g-1}\beta$, $\In(f_3^g)=\alpha^{g-1}\gamma$.
\end{prop}
\pf
By induction on $g$. $g=1$ is clear by the first line of
Proposition~\ref{recursion}. Applying our recursion relations and the
induction hypothesis, we get
$f_1^{g+1}=\alpha^{g+1}+\alpha^g\beta+\ldots\,$,
$f_2^{g+1}=\alpha^g\beta+\frac{2g}{g+1}\alpha^{g-1}\gamma+\ldots\,$,
$f_3^{g+1}=\alpha^g\gamma+\ldots\,$, where $\ldots$ mean terms of lower
order.
\qed

Now setting $\cali_g:=(f_1^g,f_2^g,f_3^g)\subset
\qz[\alpha,\beta,\gamma]$, the ideal spanned by $f_i^g$, $i=1,2,3$, then
\begin{prop}\label{inideal}
  $\In(\cali_g)=(\alpha^a\beta^b\gamma^c,a+b+c\ge g)$.
\end{prop}
\pf
By induction on $g$, $g=1$ being trivially true. From
\begin{eqnarray*}
  \gamma f_1^g&=&f_3^{g+1}\\
  g^2\gamma f_2^g&=&\gamma f_1^{g+1}-\alpha\gamma f_1^g
  \ =\ \gamma f_1^{g+1}-\alpha f_3^{g+1}\\
  \frac{4g}{g+1}\gamma f_3^g&=&\gamma f_2^{g+1}-\beta\gamma f_1^g
  \ =\ \gamma f_2^{g+1}-\beta f_3^{g+1}
\end{eqnarray*}
we see that $\gamma\cali_g\subset\cali_{g+1}$ (this is also clear from
the observation that $\gamma\in H^*(\calm_{g+1})$ is Poincar\'e dual to
$2g$ copies of $\modg$, cf.\ below). By induction hypothesis the claim
is thus true for $c>0$. But in any homogeneous expression (with
$\alpha$, $\beta$, $\gamma$ having weights $1$, $2$, $3$ respectively)
the monomials containing $\gamma$ have lower order than those without.
Therefore, we can reduce modulo $\gamma$ (indicated by a bar) and have
only to show $\In(\bar{\cali_g})=(\bar\alpha^a\bar\beta^b,a+b\ge g)$.
Modulo $\gamma$ the recursion relations read
\[
  \bar{f}_1^{g+1}=\bar\alpha\bar{f}_1^g+g^2\bar{f}_2^g,\ \ \
  \bar{f}_2^{g+1}=\bar\beta\bar{f}_1^g\,.
\]
Now we are able to repeat the argument from above with $\bar\beta$
instead of $\gamma$, because
\begin{eqnarray*}
  \bar\beta\bar{f}_1^g&=&\bar{f}_2^{g+1}\\
  g^2\bar\beta\bar{f}_2^g&=&\bar\beta\bar{f}_1^{g+1}-\bar\alpha\bar\beta
  \bar{f}_1^g\ =\ \bar\beta\bar{f}_1^{g+1}-\bar\alpha\bar{f}_2^{g+1}\,.
\end{eqnarray*}
This leaves us with the case $b=0$, $c=0$, which is clearly true since
$\alpha^g=\In(f_1^g)$ is the smallest power of $\alpha$ contained in
$\cali_g$ (for $\deg f_i^g\le g$, $i=1,2,3$).
\qed

\section{Completeness}
The strategy of showing that $\cali_g=(f_1^g,f_2^g,f_3^g)\subset
\cz[\alpha,\beta,\gamma]$ is really the ideal of relations between
$\alpha$, $\beta$, $\gamma$ is a simple dimension count.  But although
the subring $\langle\alpha,\beta,\gamma\rangle\subset H^*(\modg,\qz)$
generated by $\alpha$, $\beta$, $\gamma$ is the invariant ring under
the action of the orientation preserving diffeomorphisms, the authors
do not know a direct way to compute
$\dim_\qz\langle\alpha,\beta,\gamma\rangle$. Instead we are viewing the
even cohomology $H^{2*}(\modg,\qz)$ as module over
$\cz[\alpha,\beta,\gamma]/ \cali_g$ and check injectivity of the
structure map $\cz[\alpha,\beta,\gamma]/ \cali_g\rightarrow
H^{2*}(\modg,\qz)$ by refining the basis $\{\alpha^a
\beta^b\gamma^c\mid a+b+c<g\}$ of $\cz[\alpha,\beta,\gamma]/\cali_g$ to
a basis of $H^{2*}(\modg,\qz)$. As a by-product we will actually find
an explicit basis of the latter, which in a sense explains the
inductive formulas for the even Betti numbers found by Newstead
\cite{newstead2}.
\begin{prop}\label{generators}
  $H^{2*}(\modg,\qz)$ is generated by elements of the form\\[5pt]
  \hspace*{2cm} $\alpha^a\beta^b\psi_{i_1}\ldots\psi_{i_{2l}}$
  with $a+b+2l<g-1$, \\
  \hspace*{1cm}and $\alpha^a\beta^b\gamma^k\psi_{i_1}\ldots\psi_{i_{2l}}$
  with $a+b+k+2l=g-1$, $k\ge0$,\\[5pt]
  where $1\le i_1<\ldots<i_{2l}\le2g$.
\end{prop}
As an intermediate notion between the $\psi_i$ and $\gamma$ let us
introduce the classes $\gamma_j:=\psi_j\psi_{j+g}$, $j=1,\ldots,g$
(then $\gamma=2\sum_j\gamma_j$). Each of the $\gamma_j$ is Poincar\'e
dual to a diffeomorphic image $N_j$ of $\calm_{g-1}$ (by ``contracting
a handle'', cf.\ no.26 in \cite{thaddeus}). Moreover, $\calu|N_j$ is
topologically a universal bundle on $\calm_{g-1}$, so $\alpha$,
$\beta$, $\gamma$ restrict to generators $\hat\alpha$, $\hat\beta$,
$\hat\psi_i$ ($i\neq j,j+g$) of $H^*(\calm_{g-1},\qz)$
($\psi_j|N_j=0=\psi_{j+g}|N_j$ since $\psi_j\gamma_j=0
=\psi_{j+g}\gamma_j$ trivially). We will also use the fact that
intersection proucts $\alpha^a\beta^b\psi_{i_1}\ldots\psi_{2l}[\modg]$
($a+2b+3l=3g-3$) are zero unless
$\{i_1,\ldots,i_{2l}\}=\{j_1,j_1+g,\ldots,j_l,j_l+g\}$ in which case
\[
  \alpha^a\beta^b\psi_{i_1}\ldots\psi_{i_{2l}}[\modg]
  =\frac{1}{g!}\alpha^a\beta^b\gamma^l[\modg]\,,
\]
i.e.\ only depending on the length of the sequence $(j_1,\ldots,j_l)$,
cf.\ \cite{thaddeus}.
\vspace{15pt}

\noindent
{\em Proof of proposition.}
We want to refine the result of Proposition~\ref{inideal} that a
monomial $\alpha^a\beta^b\gamma^c$ is equivalent (= may be reduced
modulo $\cali_g$) to a polynomial of lower order. For this we will use
the reverse lexicographic order
$\alpha>\beta>\psi_1>\ldots>\psi_{2g}>\gamma_1>\ldots>\gamma_{2g}>
\gamma$.

Let $1\le i_1<\ldots<i_k\le2g$ and $1\le j_1<\ldots<j_k\le g$ with
$\{i_1,\ldots,i_k\}\cap\{j_1,j_1+g,\ldots,j_l,j_l+g\}=\emptyset$.
\vspace{5pt}

\noindent
{\em Claim:} If $a+b+k+l\ge g$ then
$\alpha^a\beta^b\psi_{i_1}\ldots\psi_{i_k}
\gamma_{j_1}\ldots\gamma_{j_l}$ is equivalent to a polynomial of lower
order modulo $\cali_g$, which can be taken of the form
$F(\alpha,\beta)\,
\psi_{i_1}\ldots\psi_{i_k}\gamma_{j_1}\ldots\gamma_{j_l}$.
\vspace{-7pt}

\noindent
The claim certainly holds if $k+l=0$ by Proposition~\ref{inideal}. For
$l>0$ let $\iota:\calm_{g-1}\hookrightarrow\modg$ have image $N_{j_k}$
(Poincar\'e dual to $\gamma_{j_k}$) and use a hat to denote pull-back
by $\iota$. By descending induction on $g$ then ($l<g$ because
$\gamma_1\ldots\gamma_g=0$ for dimension reasons),
\begin{eqnarray*}
  \lefteqn{\iota^*\left(\alpha^a\beta^b\psi_{i_1}\ldots\psi_{i_k}
  \gamma_{j_1}\ldots\gamma_{j_{l-1}}\right)
  \ =\ {\hat\alpha}^a{\hat\beta}^b{\hat\psi_{i_1}}\ldots{\hat\psi_{i_k}}
  {\hat\gamma_{j_1}}\ldots{\hat\gamma_{j_{l-1}}} }\\
  &=&F(\hat\alpha,\hat\beta)\,{\hat\psi_{i_1}}\ldots{\hat\psi_{i_k}}
  {\hat\gamma_{j_1}}\ldots{\hat\gamma_{j_{l-1}}}
  \ =\ \iota^*\left(F(\alpha,\beta)\,\psi_{i_1}\ldots\psi_{i_k}
  \gamma_{j_1}\ldots\gamma_{j_{l-1}}\right)
\end{eqnarray*}
with $\mbox{order}(F)<a+b$. This means $\alpha^a\beta^b\psi_{i_1}\ldots
\psi_{i_k}\gamma_{j_1}\ldots\gamma_{j_l}=F(\alpha,\beta)\,\psi_{i_1}
\ldots\psi_{i_k}\gamma_{j_1}$ $\ldots\gamma_{j_l}$ as wanted. Finally
the case $l=0$, $k>0$: Set $\bar\imath_k=i_k+g$ if $i_k\le g$ and
$\bar\imath_k=i_k-g$ if $i_k>g$. By the previous case, we get
$\alpha^a\beta^b\psi_{i_1}\ldots
\psi_{i_k}\psi_{\bar\imath_k}=F(\alpha,\beta)\psi_{i_1}\ldots
\psi_{i_k}\psi_{\bar\imath_k}$. Then the above remarks on the
intersection product show
\[
  \left(\alpha^a\beta^b\psi_{i_1}\ldots\psi_{i_k}
  -F(\alpha,\beta)\,\psi_{i_1}\ldots\psi_{i_k}\right)\cdot
  A[\modg]\ =\ 0
\]
for all $A\in H^*(\modg)$,
i.e.\ $\alpha^a\beta^b\psi_{i_1}\ldots\psi_{i_k}
=F(\alpha,\beta)\psi_{i_1}\ldots\psi_{i_k}$. This proves the claim.\\
The proposition is now clearly reduced to a second
\vspace{5pt}

\noindent
{\em Claim:} Let $M=\alpha^a\beta^b\psi_{i_1}\ldots\psi_{i_k}
\gamma_{j_1}\ldots\gamma_{j_l}$ with $a+b+k+l=g-2$. Then for all $1\le
i,j \le g$, $M\gamma_i-M\gamma_j$ may be reduced to lower order modulo
$\cali_g$.
\vspace{5pt}

\noindent
In fact, from the above we know already
$M\gamma_i\gamma_j=F\gamma_i\gamma_j$ in $H^*(\modg)$ with
$\mbox{order}(F)<g-2$. $F\gamma_i-F\gamma_j$ is our candidate for the
lower order term. If $\{i,i+g,j,j+g\}\cap\{i_1,\ldots,i_k\} =\emptyset$
and $A=\alpha^{a'}\beta^{b'}\psi_{i_1}\ldots\psi_{i_k}$ then
\[
  \left(M\gamma_i-M\gamma_j\right)\cdot A[\modg]\ =\ 0\
  =\ \left(F\gamma_i-F\gamma_j\right)\cdot A[\modg]
\]
again by the symmetry in the $\gamma_i$ of the intersection pairing. On
the other hand
\begin{eqnarray*}
  (M\gamma_i-M\gamma_j)\gamma_j[\modg]
  &=&M\gamma_i\gamma_j A[\modg]
  \ =\ F\gamma_i\gamma_j A[\modg]\\
  &=& (F\gamma_i-\gamma_j)\gamma_j[\modg]
\end{eqnarray*}
and analogously with $\gamma_i$. Thus $M\gamma_i=M\gamma_j+
F(\gamma_i-\gamma_j)$ modulo $\cali_g$ as claimed.\\
(Note: This argument fails for the $\psi_i$ because of the presence of
$\psi_{i+g}$ respectively $\psi_{i-g}$.)
\qed

The only thing we finally need to do is to count the number of generators
and compare with the inductive formula for the Betti numbers found by
Newstead.
\begin{prop}\label{dimcount}
  The generators for $H^{2*}(\modg,\qz)$ in Proposition~\ref{generators}
  are linearly independent up to the middle dimension.
\end{prop}
\pf
Let $G_r$ be the set of generators of (real) degree $r$ from
Proposition~\ref{generators}, $g_r:=\sharp G_r$. We will show that for
$s\displaystyle\le\left[\frac{3g-1}{2}\right]$
\[
  g_{2s+4}\ =\ g_{2s}+\sum_{l=s-g+1}^{[s/3]}{2g\choose 2l}
\]
which together with $g_0=1$ and $g_2=1$ is exactly Newstead's formula
for the even Betti numbers \cite{newstead1}.
Define a map $\ph:G_{2s}\rightarrow G_{2s+4}$ by
\[
  \alpha^a\beta^b\psi_{i_1}\ldots\psi_{i_{2l}}\longmapsto
  \alpha^a\beta^{b+1}\psi_{i_1}\ldots\psi_{i_{2l}}\ \ \
  \mbox{for }a+b+2l<g-1
\]
and
\begin{eqnarray*}
  \alpha^a\beta^b\gamma^k\psi_{i_1}\ldots\psi_{i_{2l}}\longmapsto
  \alpha^{a-1}\beta^b\gamma^{k+1}\psi_{i_1}\ldots\psi_{i_{2l}}
  &&\mbox{ for }a+b+k+2l=g-1,
\end{eqnarray*}
with $k\ge0$ and $a>0$.
Note that the case $a=k=0$ does not occur (then $b+2l=g-1$ and
$2b+3l=s$ imply $s-2g+2=-l\le0$ which contradicts
$s\le(3g-1)/2$). Now $G_{2s+4}\setminus\img\ph=\{\alpha^a\psi_{i_1}\ldots
\psi_{i_{2l}}\mid a+2l\le g-1, a+3l=s\}$ s.th.\ $l$ runs from
$s-g+1$ to $[s/3]$ ($a$ is determined through $a+3l=s$). The contribution
for $l$ fixed is then precisely $2g\choose2l$.
\qed

\end{document}